\documentclass{article}

\usepackage{arxiv}

\usepackage{color}
\usepackage[colorlinks,linkcolor=blue,citecolor=blue]{hyperref}
\usepackage{amsmath,amsfonts,amssymb,amsthm,mathrsfs}
\usepackage[margin=10pt,font=small,labelfont=bf,labelsep=period]{caption}

\usepackage{indentfirst}  
\usepackage{cite}

\usepackage[utf8]{inputenc} 
\usepackage[T1]{fontenc}    
\usepackage{url}            
\usepackage{booktabs}       
\usepackage{nicefrac}       
\usepackage{graphicx,subfigure}
\usepackage{doi}

\usepackage{algorithm}
\usepackage{listings,xcolor} 
\usepackage{courier}
\usepackage{textcomp}
\usepackage[framed,numbered,autolinebreaks,useliterate]{mcode} 

\DeclareCaptionFont{white}{\color{white}}
\DeclareCaptionFormat{listing}{\colorbox[cmyk]{0.43, 0.35, 0.35,0.01}{\parbox{\textwidth}{#1#2#3}}}
\title{Implementation of Polygonal Mesh Refinement in MATLAB}

\author{ Yue Yu \\
	\texttt{terenceyuyue@sjtu.edu.cn} \\
}

\date{}

\renewcommand{\headeright}{}
\renewcommand{\undertitle}{}
\renewcommand{\shorttitle}{Polygonal Mesh Refinement}

\begin{document}
\maketitle

\begin{abstract}
  We present a simple and efficient MATLAB implementation of the local refinement for polygonal meshes. The purpose of this implementation is primarily educational, especially on the teaching of adaptive virtual element methods.
\end{abstract}


\section{Introduction}
The virtual element method (VEM), first introduced in \cite{Beirao-Brezzi-Cangiani-2013}, is a generalization
of the standard finite element method on general polytopal meshes. In the past few years, people have witnessed
rapid progress of virtual element methods (VEMs) for numerically solving various partial differential equations, see for
examples \cite{Ahmad-Alsaedi-Brezzi-2013,DeDios-Lipnikov-Manzini-2016,Zhao-Zhang-Chen-2018,Beirao-DaVeiga-Brezzi-2014}.
One can also refer to \cite{Sutton-2017} for a transparent MATLAB implementation of the conforming
virtual element method for Poisson equation. Due to the large flexibility of the meshes, researchers gradually turn their attention to the posterior error analysis of the VEMs and have made some progress in recent years
(cf. \cite{Beirao-Manzini-2015,Cangiani-Georgoulis-Pryer-2017,Berrone-Borio-2017,Chi-Beirao-Paulino-2019,Beirao-Manzini-Mascotto-2019}).

As we all know, standard adaptive algorithms based on the local mesh refinement can be written as the following loops
\[{\bf SOLVE} \to {\bf ESTIMATE} \to {\bf MARK} \to {\bf REFINE}.\]
Given an initial polygonal subdivision $\mathcal{T}_0$, to get $\mathcal{T}_{k+1}$ from $\mathcal{T}_k$ we
first solve the VEM problem under consideration to get the numerical solution $u_k$ on $\mathcal{T}_k$. The error is then estimated by using $u_k$, $\mathcal{T}_k$ and the a posteriori error bound. And the local error bound is used to mark a subset
of elements in $\mathcal{T}_k$ for refinement. The marked polygons and possible more neighboring elements are refined
in such a way that the subdivision meets certain conditions, such as shape regularity. In the implementation, it is usually time-consuming to write a mesh refinement routine since we need to carefully design the rule for dividing the marked elements to get a refined mesh of high quality.

In this paper, we are intended to present an efficient MATLAB implementation of the mesh refinement for polygonal meshes.
To the best of our knowledge, this is the first publicly available implementation of the polygonal mesh refinement algorithms.
We divide elements by connecting the midpoint of each edge to its barycenter, which may be the most natural partition frequently used in VEM papers, referred to as 4-node subdivision in this context.
To remove small edges, some additional neighboring polygons of the marked elements are included in the refinement set by requiring the one-hanging-node rule: limit the mesh to have just one hanging node per edge.
We discuss the implementation step by step and give an application in the posteriori error analysis
for Poisson equation in the last section. The current implementation or the 4-node subdivision requires
that the barycenter is an internal point of each element.

\section{Strategy for removing small edges} \label{sect:remove}

\begin{figure}[!htb]
  \centering
  \subfigure[]{\includegraphics[scale=0.5]{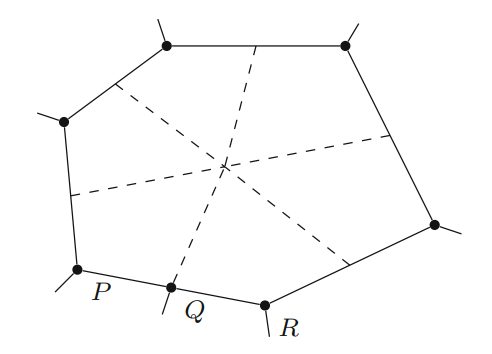}} \hspace{1cm}
  \subfigure[]{\includegraphics[scale=0.45,trim=0 90 0 80,clip]{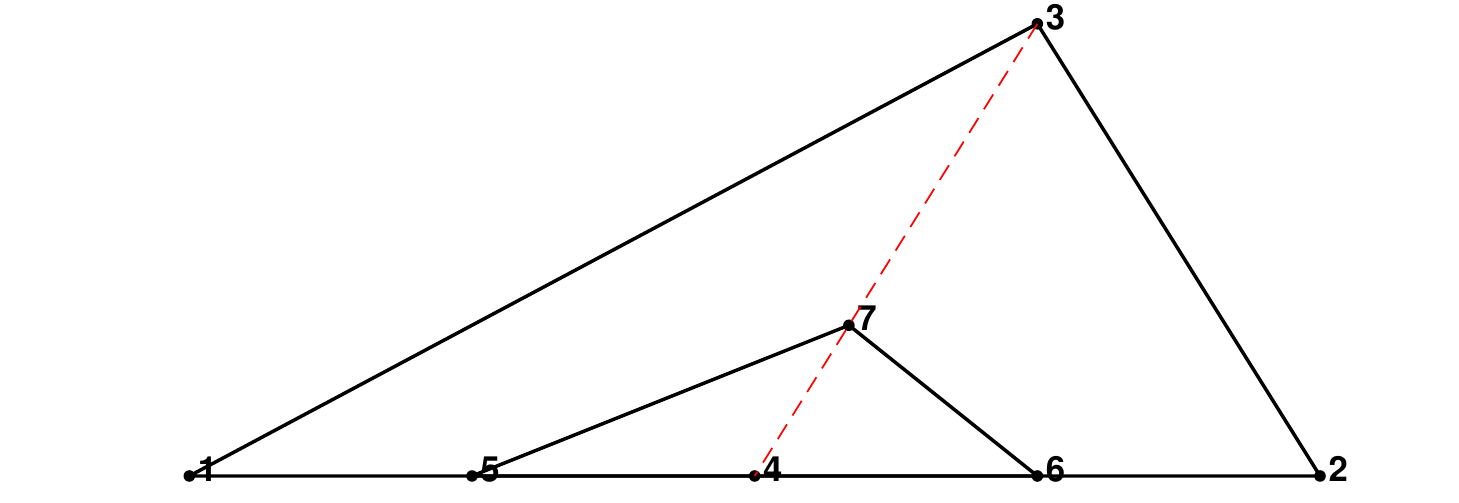}}\\
 \caption{Edge bisection. (a) Admissible bisection: add a new edge from hanging node to the barycenter;
 (b) Two-edge bisection: cut two edges with hanging node as endpoint}\label{fig:refine0}
\end{figure}

For polygonal meshes, a natural mesh refinement strategy may be the 4-node subdivision as illustrated in
Fig.~\ref{fig:refine0}~(a) given in \cite{Cangiani-Georgoulis-Pryer-2017}, where a target element is divided
by connecting the midpoint of each edge to the barycenter. For the element with hanging node $Q$, it is better
to add a new edge from $Q$ to the barycenter instead of bisecting two edges $PQ$ and $QR$, otherwise degenerate
triangles may appear. We refer to this treatment as the admissible bisection. As shown in Fig.~\ref{fig:refine0}~(b),
the triangle $\Delta 123$ is generated by bisecting $PQ$ and $QR$ and the new triangle $\Delta 567$ is obtained
from $\Delta 123$ in a similar manner. According to the properties of triangles, the barycenter $z_7$ lies on
the median line $e_{34}$ and $|e_{47}|:|e_{73}| = 1:2$, which implies that $\angle 7 > \angle 3 $ and hence a
degenerate triangle will appear if we continue with the ``two-edge bisection''. On the other hand,
the admissible bisection leads to a simpler polygonal mesh.

A mesh refinement function has the following form

\vspace{0.5em}

  \hspace{2cm} \mcode{[node,elem] = PolyMeshRefine(node,elem,idElemMarked)}.

\vspace{0.5em}

\noindent Here, \mcode{node} and \mcode{elem} are two basic data structure, storing the coordinates of nodes and the
connectivity list, respectively, and the array \mcode{idElemMarked} gives the index of marked elements.
Generally, the elements to be refined are more than marked elements. By imposing no restriction on the number
of hanging nodes on a single edge, we are at risk of violating the ``no short edge'' assumption
(cf. \cite{Cangiani-Georgoulis-Pryer-2017}). Although this requirement does not seem to be necessary for the VEMs to
remain accurate and stable in practice, it is still expected to produce high quality meshes without small edges. To this end,
some additional cells should be included in the refinement set.

\begin{figure}[!htb]
  \centering
  \subfigure[]{\includegraphics[scale=0.45,trim=0 50 0 50,clip]{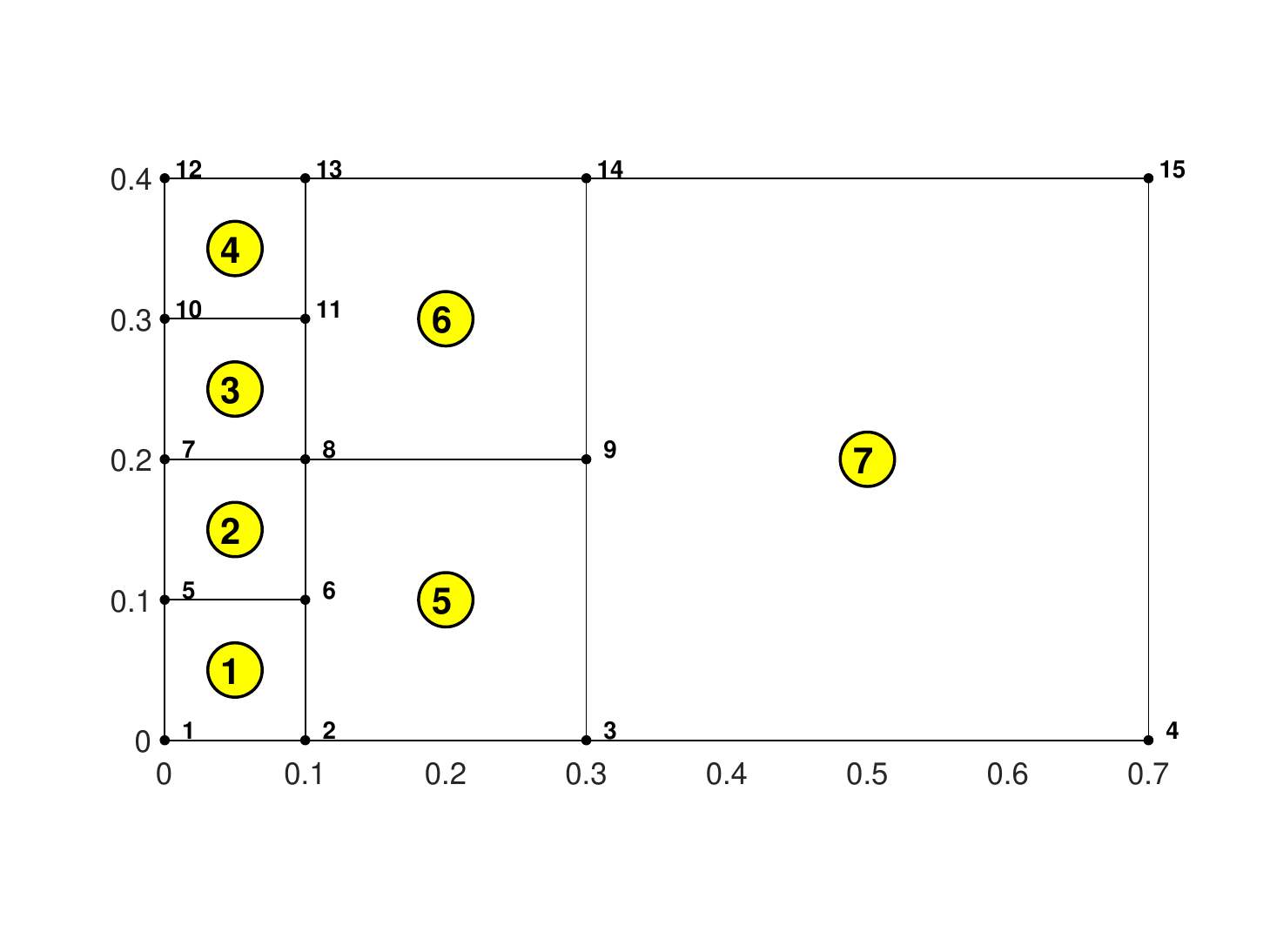}}
  \subfigure[]{\includegraphics[scale=0.45,trim=0 50 0 50,clip]{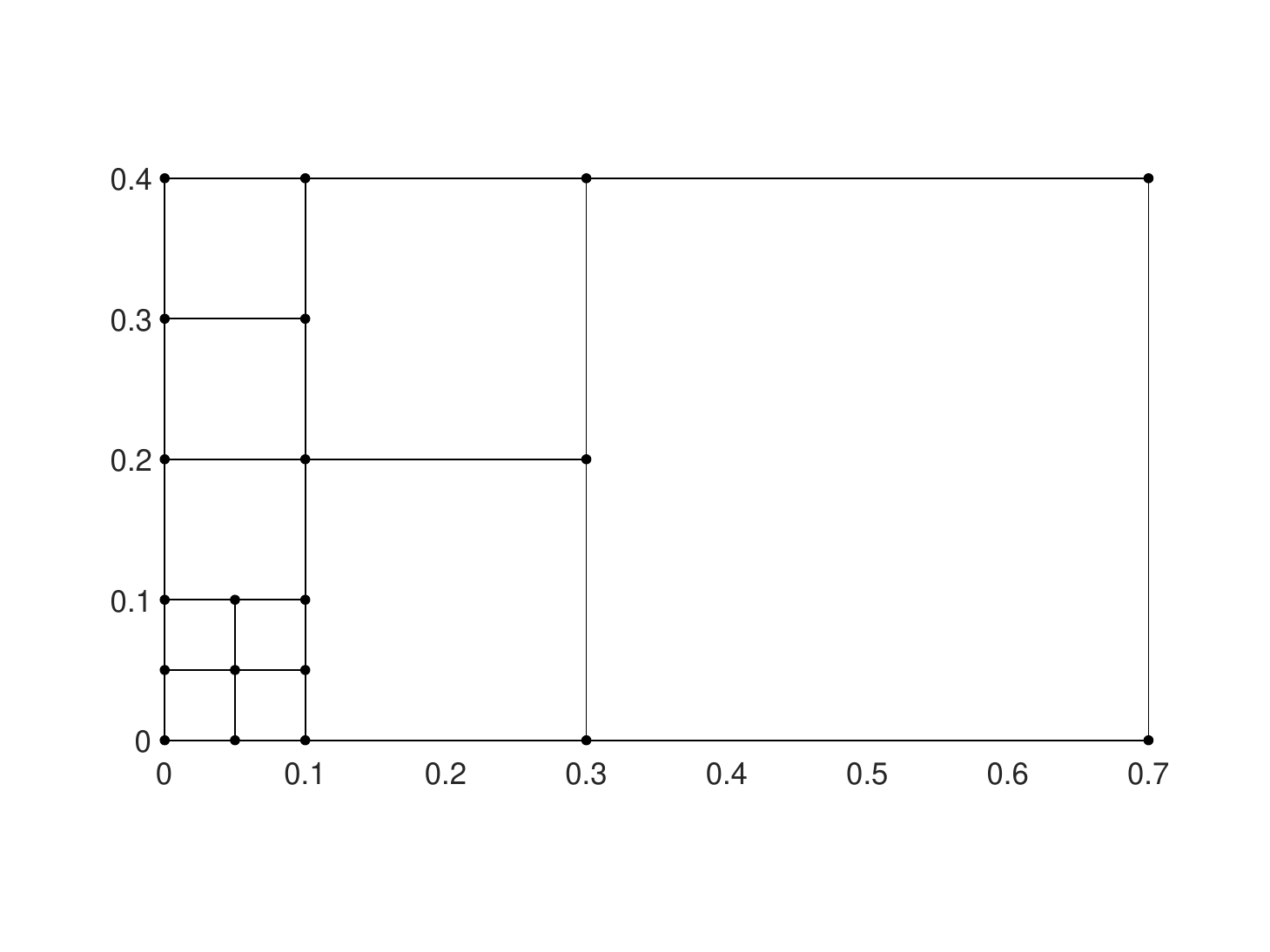}}\\
\caption{(a) Initial mesh; (b) Refinement of only the marked elements.}\label{fig:refine01}
\end{figure}

\begin{figure}[!htb]
  \centering
  \includegraphics[scale=0.45,trim=0 50 0 50,clip]{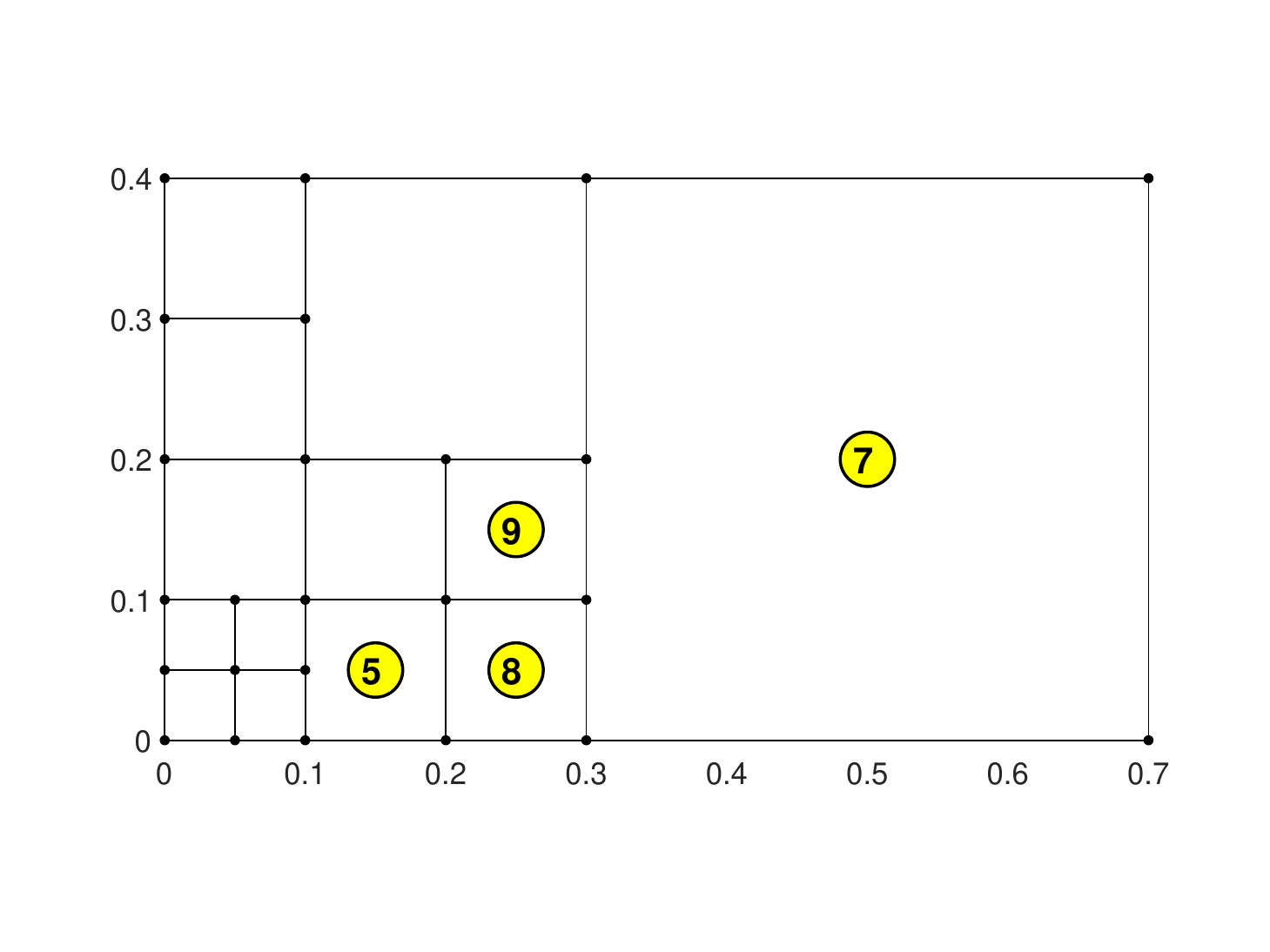}\\
\caption{Refine both \textcircled{1} and \textcircled{5}}\label{fig:refineh2}
\end{figure}

Given an initial mesh shown in Fig.~\ref{fig:refine01}~(a), we find that the short edge problem reduces to the refinement
of the two adjacent polygons \textcircled{1} and  \textcircled{5} as a typical case.
\begin{itemize}
  \item Suppose that \textcircled{1} is a marked element: \mcode{idElemMarked = 1}. If we only refine the marked element
  as shown in Fig.~\ref{fig:refine01}~(b) and proceed to divide the small cells in the lower right corner, then small
  edges appear in \textcircled{5} due to the frequently added hanging nodes. This phenomenon also happens in the two-edge
  bisection illustrated in Fig.~\ref{fig:refine0}~(b), which is resolved by using the admissible bisection.
  \item To avoid the occurrence of short edges, we further refine the adjacent polygon \textcircled{5}, yielding a new mesh given
  in Fig.~\ref{fig:refineh2}. That is, for two elements $K_1$ and $K_2$ sharing an edge $e$, we refine both $K_1$ and $K_2$
  if $K_1$ is in the refinement set and one endpoint of $e$ is the hanging node of $K_2$.
  \item For elements in Fig.~\ref{fig:refine01}~(a), \textcircled{5} and \textcircled{7} can be viewed as \textcircled{1} and
  \textcircled{5}, respectively. For the same reason, we also need to partition the polygon \textcircled{7}.
  \item Repeating the above procedure, one can collect all the new elements for refinement and hence avoid producing small
  edges since the above treatment ensures that each edge in the resulting mesh contains at most one hanging node
  (a midpoint of the collinear edge).
\end{itemize}

In the following, an edge is called nontrivial if one of its endpoints is hanging node.

\section{Implementation}

\subsection{Data structure}

We first discuss the data structure to represent polygonal meshes so as to facilitate the refinement procedure. The idea stems from the treatment of triangulation in iFEM (cf. \cite{iFEM}), which is generalized to polygonal meshes with certain modifications.

We will only maintain and update two basic data structure \mcode{node} and \mcode{elem}. In the node matrix \mcode{node}, the first and second rows contain $x$- and $y$-coordinates of the nodes in the mesh. The cell array \mcode{elem} records the vertices of each element with a counterclockwise order.

We always use the symbols \mcode{N}, \mcode{NT}, and \mcode{NE} to represent the number of nodes, triangles, and edges, respectively. The following code will build auxiliary data structure.
\vspace{-0.5cm}
\begin{lstlisting}
%% Get auxiliary data
NT = size(elem,1);
if ~iscell(elem), elem = mat2cell(elem,ones(NT,1),length(elem(1,:))); end
% centroid
centroid = zeros(NT,2); diameter = zeros(NT,1);
s = 1;
for iel = 1:NT
    index = elem{iel};
    verts = node(index,:); verts1 = verts([2:end,1],:);
    area_components = verts(:,1).*verts1(:,2)-verts1(:,1).*verts(:,2);
    ar = 0.5*abs(sum(area_components));
    centroid(s,:) = sum((verts+verts1).*repmat(area_components,1,2))/(6*ar);
    diameter(s) = max(pdist(verts));
    s = s+1;
end
if max(diameter)<4*eps, error('The mesh is too dense'); end
% totalEdge
shiftfun = @(verts) [verts(2:end),verts(1)];
T1 = cellfun(shiftfun, elem, 'UniformOutput', false);
v0 = horzcat(elem{:})'; v1 = horzcat(T1{:})';
totalEdge = sort([v0,v1],2);
% edge, elem2edge
[edge, i1, totalJ] = unique(totalEdge,'rows');
elemLen = cellfun('length',elem);
elem2edge = mat2cell(totalJ',1,elemLen)';
% edge2elem
Num = num2cell((1:NT)');    Len = num2cell(elemLen);
totalJelem = cellfun(@(n1,n2) n1*ones(n2,1), Num, Len, 'UniformOutput', false);
totalJelem = vertcat(totalJelem{:});
[~, i2] = unique(totalJ(end:-1:1),'rows');
i2 = length(totalEdge)+1-i2;
edge2elem = totalJelem([i1,i2]);
% neighbor
neighbor = cell(NT,1);
for iel = 1:NT
    index = elem2edge{iel};
    ia = edge2elem(index,1); ib = edge2elem(index,2);
    ia(ia==iel) = ib(ia==iel);
    neighbor{iel} = ia';
end
% number
N = size(node,1); NE = size(edge,1);
\end{lstlisting}

In the edge matrix \mcode{edge}, the first and second rows contain indices of
the starting and ending points. The row of \mcode{edge} is sorted in such a way that \mcode{edge(k,1)<edge(k,2)}
for the $k$-th edge. The cell array \mcode{elem2edge} establishes
the map of local index of edges in each polygon to its global index in matrix \mcode{edge}.
We also use the cell array \mcode{neighbor} to record the neighboring polygons for each element.
The coordinates of barycenters of all elements are stored in matrix \mcode{centroid}.

\begin{algorithm}[!htb]
\caption{Find additional elements for refinement} \label{alg:findadd}
 \begin{itemize}
  \item[-] Initialize \mcode{idElemMarkedNew} and \mcode{idElemNew} as the set of marked elements \mcode{idElemMarked}.
  \item[-] Find the adjacent elements, stored in \mcode{idElemNewAdj}, of new added elements in \mcode{idElemNew} and
  obtain the numbers of edges of all elements in \mcode{idEdgeMarkedNew}.
  \item[-] Update \mcode{idElemNew} by looping for elements in \mcode{idElemNewAdj}: Determine the index numbers of nontrivial edges in current
  element, denoted by \mcode{idEdgeDg}. If \mcode{idEdgeDg} intersects with \mcode{idEdgeMarkedNew},
  then the current element needs to be refined.
  \item[-] Update \mcode{idElemMarkedNew} by adding the elements in \mcode{idElemNew}.
  \item[-] If \mcode{idElemNew} is empty, stop; otherwise, go back to the first step.
  All the additional elements are then given by \mcode{idElemAdjRefine = setdiff(idElemMarkedNew,idElemMarked)}.
\end{itemize}
\end{algorithm}

\subsection{The additional elements for refinement}

We now collect the additional elements to be refined. Let \mcode{idElemMarked} record the marked elements
and \mcode{idElemNew} record the newly added elements in each step.
The vector \mcode{idElemMarkedNew} collects the marked elements and all the additional elements up to the current step.
The discussion in Sect. \ref{sect:remove} is summarized in Algorithm \ref{alg:findadd}.

We list the code below and explain it briefly.
\vspace{-0.5cm}
\begin{lstlisting}
%% Find the additional elements to be refined
% initialized as marked elements
idElemMarkedNew = idElemMarked; % marked and all new elements
idElemNew = idElemMarked; % new elements
while ~isempty(idElemNew)
    % adjacent polygons of new elements
    idElemNewAdj = unique(horzcat(neighbor{idElemNew}));
    idElemNewAdj = setdiff(idElemNewAdj,idElemMarkedNew);
    % edge set of new marked elements
    idEdgeMarkedNew = unique(horzcat(elem2edge{idElemMarkedNew}));
    % find the adjacent elements to be refined
    nElemNewAdj = length(idElemNewAdj);
    isRefine = false(nElemNewAdj,1);
    for s = 1:nElemNewAdj
        % current element
        iel = idElemNewAdj(s);
        index = elem{iel};  indexEdge = elem2edge{iel}; Nv = length(index);
        % local logical index of elements with hanging nodes
        v1 = [Nv,1:Nv-1]; v0 = 1:Nv; v2 = [2:Nv,1]; % left,current,right
        p1 = node(index(v1),:); p0 = node(index(v0),:); p2 = node(index(v2),:);
        err = sqrt(sum((p0-0.5*(p1+p2)).^2,2));
        ism = (err<eps);  % is midpoint
        % start the next loop if no hanging nodes exist
        if sum(ism)<1, continue; end
        % index numbers of edges connecting hanging nodes in the
        % adjacent elements to be refined
        idEdgeDg = unique(indexEdge([v1(ism),v0(ism)]));
        % whether or not the above edges are in the edge set of new marked elements
        if intersect(idEdgeDg, idEdgeMarkedNew), isRefine(s) = true; end
    end
    idElemNew = idElemNewAdj(isRefine);
    idElemMarkedNew = unique([idElemMarkedNew(:); idElemNew(:)]);
end
idElemAdjRefine = setdiff(idElemMarkedNew,idElemMarked);
\end{lstlisting}

The strategy in Sect. \ref{sect:remove} ensures that each edge has at most one hanging node for an initial mesh of high quality and the hanging node coincides with the midpoint of the collinear edge. For this reason, we can find the hanging node in a given element by computing the following errors
\[{\rm err}(i) = \Big| z_i - \frac{1}{2}(z_{i-1}+z_{i+1}) \Big|, \quad i = 1,\cdots,N_v,\]
where $z_i$ are the vertices and $N_v$ is the vertex number, as coded in Line 19-22. We remark that this process will also be used in the element refinement and element extension introduced in the following.

\subsection{Refinement of the additional elements}

We observe from Fig.~\ref{fig:refineh2} that the hanging node will appear in the subcells of the partitioned additional elements.
For this reason, we first divide the additional elements and then extend all possible elements together by adding hanging nodes.
The other reason is that we need the data structure \mcode{elem2edge} in the element extension since the midpoints
will be labeled by using the edge index.

Note that some midpoints of edges and barycenters need to be added to the matrix \mcode{node}.
We relabel the vertices, edges and elements in the following order
\[z_1,\cdots,z_{\rm N};~~e_1,\cdots,e_{\rm NE};~~K_1,\cdots,K_{\rm NT}\]
with a single index $i=1,2,\cdots,{\rm N}+{\rm NE} + {\rm NT}$, referred to as the connection number.
However, in most cases it is more convenient to use the index number in matrix \mcode{edge}.

To construct the 4-node subcells, first consider an example with the connection numbers of vertices and edges listed as
\[z_1, e_1, {\color{red} z_2}, e_2, z_3, e_3, z_4, e_4, {\color{red} z_5}, e_5, \quad N_v = 5, \]
where $z_2$ and $z_5$ are hanging nodes and the subscript stands for local index.
Next, we replace the connection numbers of nontrivial edges by the ones of hanging nodes as
\begin{align*}
                 & z_1, e_1, {\color{red} z_2}, e_2, z_3, e_3, z_4, e_4, {\color{red} z_5}, e_5 \\
 \hookrightarrow~ & z_1, {\color{red} z_2, z_2, z_2}, z_3, e_3, z_4, {\color{red} z_5, z_5, z_5}, \quad N_v = 5.
  \end{align*}
Then we can construct the corresponding data structure \mcode{elem} in a unified way, which applies to the marked elements with or
without hanging nodes.

For the data structure \mcode{elem2edge}, we simply record or modify the index numbers of trivial edges by zero.

The code is listed as follows.
\vspace{-0.5cm}
\begin{lstlisting}
%% Partition the adjacent elements to be refined
nAdjRefine = length(idElemAdjRefine);
elemAdjRefine = cell(nAdjRefine,1);
elem2edgeAdjRefine = cell(nAdjRefine,1);
for s = 1:nAdjRefine
    % current element
    iel = idElemAdjRefine(s);
    index = elem{iel}; indexEdge = elem2edge{iel}; Nv = length(index);
    % find midpoint
    v1 = [Nv,1:Nv-1]; v0 = 1:Nv; v2 = [2:Nv,1];
    p1 = node(index(v1),:); p0 = node(index(v0),:); p2 = node(index(v2),:);
    err = sqrt(sum((p0-0.5*(p1+p2)).^2,2));
    ism = (err<eps);
    % modify the edge number
    ide = indexEdge+N;  % the connection number
    ide(v1(ism)) = index(ism); ide(ism) = index(ism);
    % elem
    nsub = Nv-sum(ism);
    z1 = ide(v1(~ism));  z0 = index(~ism);
    z2 = ide(~ism);      zc = iel*ones(nsub,1)+N+NE;
    elemAdjRefine{s} = [z1(:), z0(:), z2(:), zc(:)];
    % elem2edge
    ise = false(Nv,1); ise([v1(ism),v0(ism)]) = true;
    idg = zeros(Nv,1); idg(ise) = indexEdge(ise);
    e1 = idg(v1(~ism));     e2 = idg(~ism);
    e3 = zeros(nsub,1);     e4 = zeros(nsub,1);  % numbered as 0
    elem2edgeAdjRefine{s} = [e1(:), e2(:), e3(:), e4(:)];
end
addElem = vertcat(elemAdjRefine{:});
addElem2edge = vertcat(elem2edgeAdjRefine{:}); % transform to cell arrays
if ~isempty(addElem) % may be empty
    addElem = mat2cell(addElem,ones(size(addElem,1),1),4);
    addElem2edge = mat2cell(addElem2edge,ones(size(addElem2edge,1),1),4);
end
\end{lstlisting}

\subsection{Element extension by adding hanging nodes}

The elements for extension are composed of some neighboring elements of the ones in refinement set and
some subcells of additional elements to be refined. Denote by \mcode{idEdgeCut} the index numbers of all trivial edges in
elements for refinement. Then at least one edge of the extension element corresponds to the index number in \mcode{idEdgeCut}.
We first derive the vector \mcode{idEdgeCut} as follows.
\vspace{-0.5cm}
\begin{lstlisting}
%% Extend elements by adding hanging nodes
% elements to be refined
idElemRefine = [idElemAdjRefine(:); idElemMarked(:)]; % the order cannot be changed
nRefine = length(idElemRefine);
% natural numbers of edges without hanging nodes
isEdgeCut = false(NE,1);
for s = 1:nRefine
    iel = idElemRefine(s);
    index = elem{iel}; indexEdge = elem2edge{iel}; Nv = length(index);
    v1 = [Nv,1:Nv-1]; v0 = 1:Nv; v2 = [2:Nv,1];
    p1 = node(index(v1),:); p0 = node(index(v0),:); p2 = node(index(v2),:);
    err = sqrt(sum((p0-0.5*(p1+p2)).^2,2));
    ism = (err<eps);
    idx = true(Nv,1);
    idx(v1(ism)) = false; idx(ism) = false;
    isEdgeCut(indexEdge(idx)) = true;
end
idEdgeCut = find(isEdgeCut);
\end{lstlisting}

To extend an element, we first generate a zero vector of length $2N_v$. In the odd place, we insert the vertex numbers,
while only insert index numbers of edges in \mcode{idEdgeCut} in the even place. We then obtain the connectivity vector
by deleting the zero entries.
\vspace{-0.5cm}
\begin{lstlisting}
% adjacent polygons of elements to be refined
idElemRefineAdj = unique(horzcat(neighbor{idElemRefine}));
idElemRefineAdj = setdiff(idElemRefineAdj,idElemRefine);
% basic data structure of elements to be extended
elemExtend = [elem(idElemRefineAdj); addElem];
elem2edgeExtend = [elem2edge(idElemRefineAdj); addElem2edge];
% extend the elements
for s = 1:length(elemExtend)
    index = elemExtend{s}; indexEdge = elem2edgeExtend{s};
    Nv = length(index);
    [idm,id] = intersect(indexEdge,idEdgeCut);
    idvec = zeros(1,2*Nv);
    idvec(1:2:end) = index;  idvec(2*id) = idm+N;
    elemExtend{s} = idvec(idvec>0);
end
% replace the old elements
nRefineAdj = length(idElemRefineAdj);
elem(idElemRefineAdj) = elemExtend(1:nRefineAdj);
addElem = elemExtend(nRefineAdj+1:end);
elem(idElemAdjRefine) = addElem(1:nAdjRefine);
addElem = addElem(nAdjRefine+1:end);
\end{lstlisting}

It should be pointed out that in Line 5 all the neighboring elements and all the subcells are grouped into the
vector \mcode{elemExtend}. The redundant elements do not affect the result since we simply record or modify the index numbers of
trivial edges of additional elements by zero.

\subsection{Partition of the marked elements}

We can refine the marked elements in the same way as the additional elements for refinement.
Note that in obtaining the vector \mcode{idEdgeCut}, the original data structure \mcode{elem} of marked elements are needed.
Therefore, it is necessary to partition the marked elements after element extension.
\vspace{-0.5cm}
\begin{lstlisting}
%% Partition the marked elements
nMarked = length(idElemMarked);
addElemMarked = cell(nMarked,1);
for s = 1:nMarked
    % current element
    iel = idElemMarked(s);
    index = elem{iel}; indexEdge = elem2edge{iel}; Nv = length(index);
    % find midpoint
    v1 = [Nv,1:Nv-1]; v0 = 1:Nv; v2 = [2:Nv,1];
    p1 = node(index(v1),:); p0 = node(index(v0),:); p2 = node(index(v2),:);
    err = sqrt(sum((p0-0.5*(p1+p2)).^2,2));
    ism = (err<eps);
    % replace the edge numbers with the numbers of hanging nodes
    ide = indexEdge+N;  % connection number
    ide(v1(ism)) = index(ism); ide(ism) = index(ism);
    % partition the elements with or without hanging nodes
    nsub = Nv-sum(ism);
    z1 = ide(v1(~ism));   z0 = index(~ism);
    z2 = ide(~ism);       zc = iel*ones(nsub,1)+N+NE;
    addElemMarked{s} = [z1(:), z0(:), z2(:), zc(:)];
end
% replace the old elements
addElemMarked = vertcat(addElemMarked{:});
addElemMarked = mat2cell(addElemMarked, ones(size(addElemMarked,1),1), 4);
elem(idElemMarked) = addElemMarked(1:nMarked);
addElemMarked = addElemMarked(nMarked+1:end);
\end{lstlisting}

\subsection{Update of the basic data structure}

We finally update the basic data structure \mcode{node} and \mcode{elem} by adding all new elements and reordering the vertices.

\vspace{-0.5cm}
\begin{lstlisting}
%% Update node and elem
idElemRefine = unique(idElemRefine); % in ascending order
z1 = node(edge(idEdgeCut,1),:); z2 = node(edge(idEdgeCut,2),:);
nodeEdgeCut = (z1+z2)/2;
nodeCenter = centroid(idElemRefine,:);
node = [node; nodeEdgeCut; nodeCenter];
elem = [elem; addElem; addElemMarked];
%% Reorder the vertices
[~,~,totalid] = unique(horzcat(elem{:})');
elemLen = cellfun('length',elem);
elem = mat2cell(totalid', 1, elemLen)';
\end{lstlisting}

\section{Application in the posteriori error estimates for virtual element
methods}

The mesh refinement routine \mcode{PolyMeshRefine.m} is formed by collecting all the code fragments in the last section, which
is also available from GitHub (\url{https://github.com/Terenceyuyue/mVEM}) as part of the mVEM package containing efficient and
easy-following codes for various VEMs published in the literature.
The refinement function has been tested for many initial meshes generated by PolyMesher (cf. \cite{PolyMesher-2012}), a polygonal mesh
generator based on the Centroidal Voronoi Tessellations (CVTs), which confirms the effectiveness and correctness of our code.
We omit the details and in turn consider the application in the posteriori error estimates of the VEM for Poisson equation.

Consider the Poisson equation with Dirichlet boundary condition on the unit square. The exact solution is given by
\[u(x,y) = xy(1 - x)(1 - y){\text{exp}}\left(  - 1000((x - 0.5)^2 + (y - 0.117)^2) \right).\]
The error estimator is taken from \cite{Berrone-Borio-2017}. We employ the VEM in the lowest order case and use the
D\"{o}rfler marking strategy with parameter $\theta = 0.4$ to select the subset of elements for refinement.

\begin{figure}[!htb]
  \centering
  \subfigure[]{\includegraphics[scale=0.45]{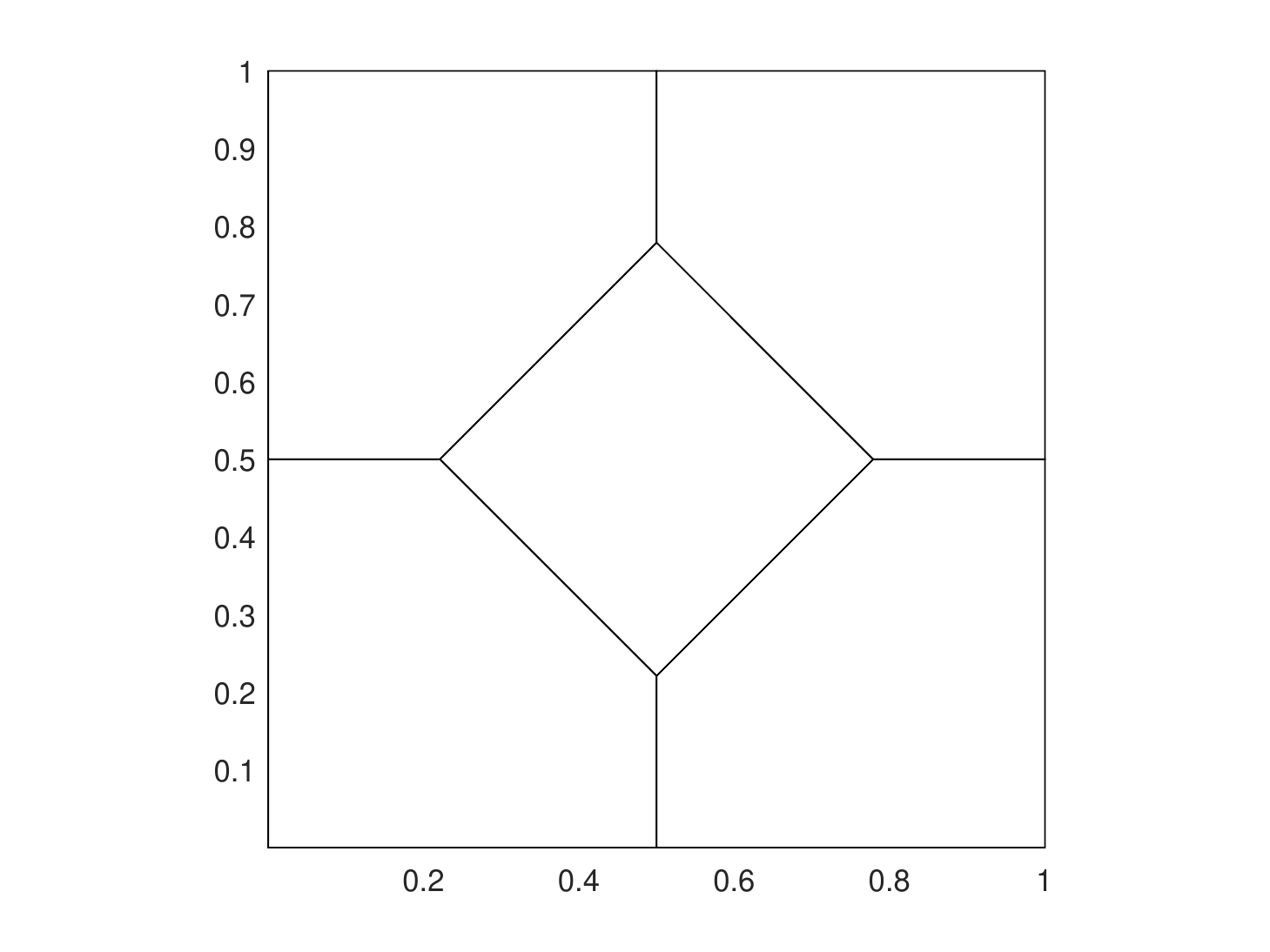}}
  \subfigure[]{\includegraphics[scale=0.45]{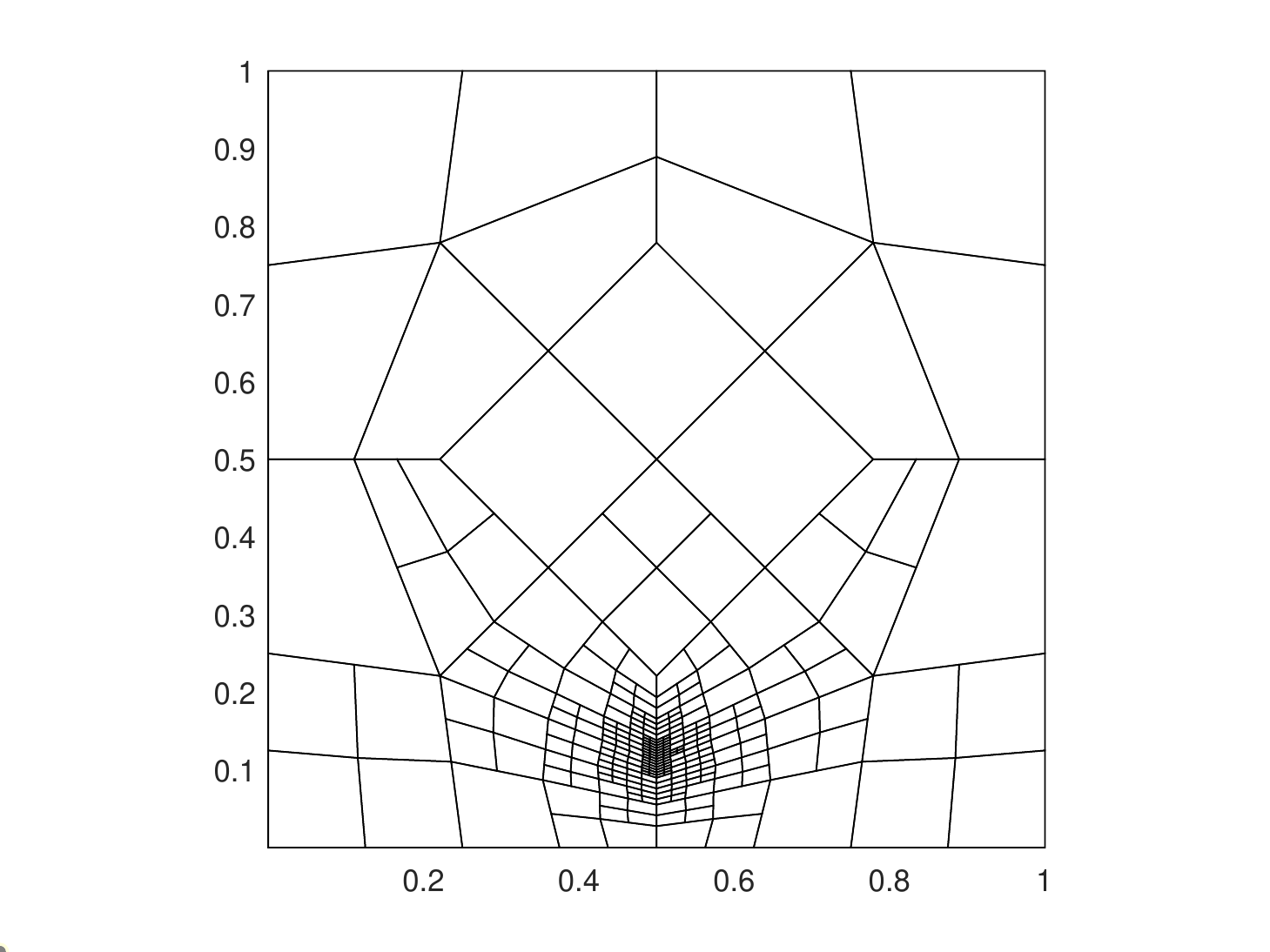}}\\
  \subfigure[]{\includegraphics[scale=0.45]{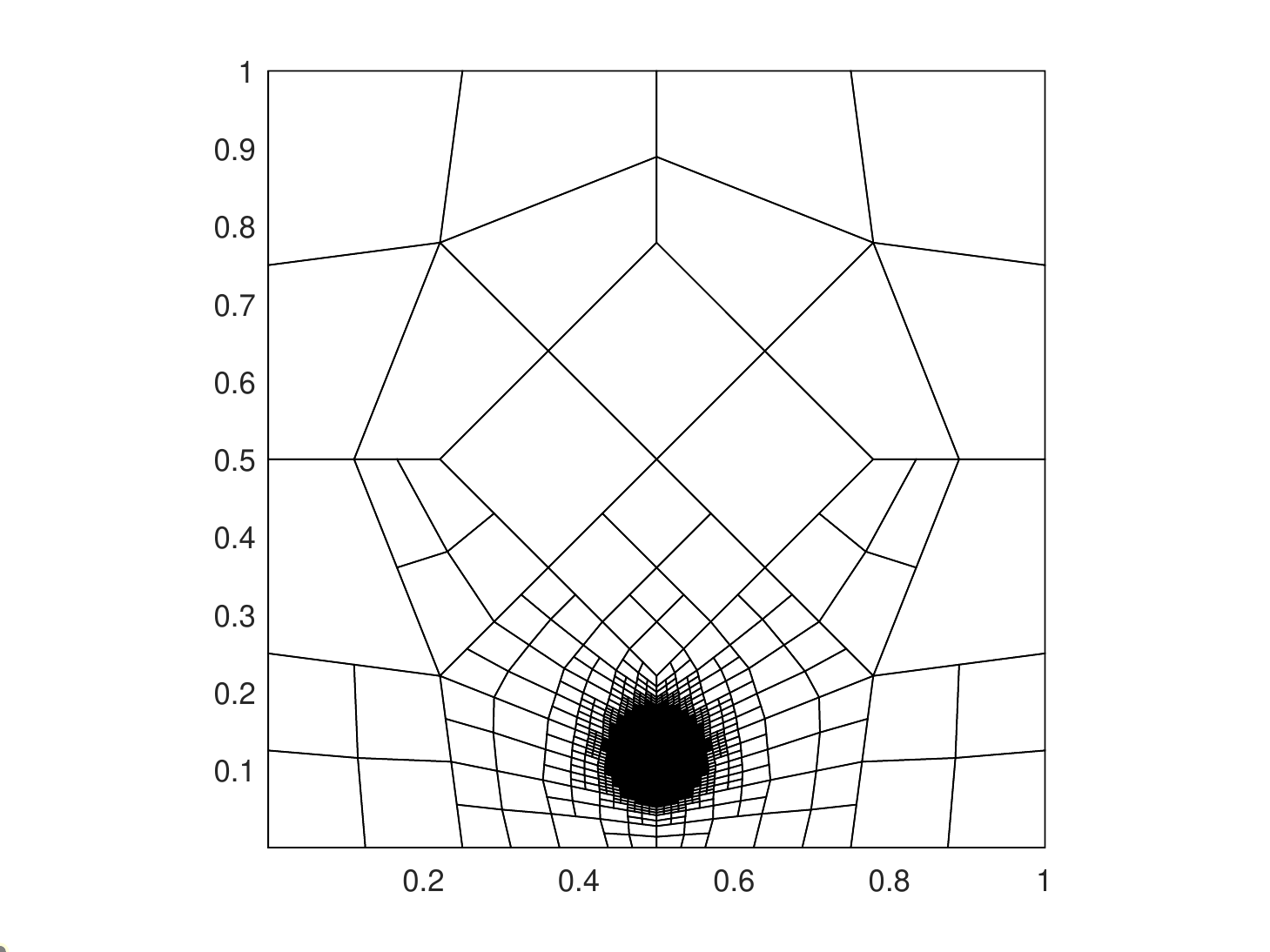}}
  \subfigure[]{\includegraphics[scale=0.45]{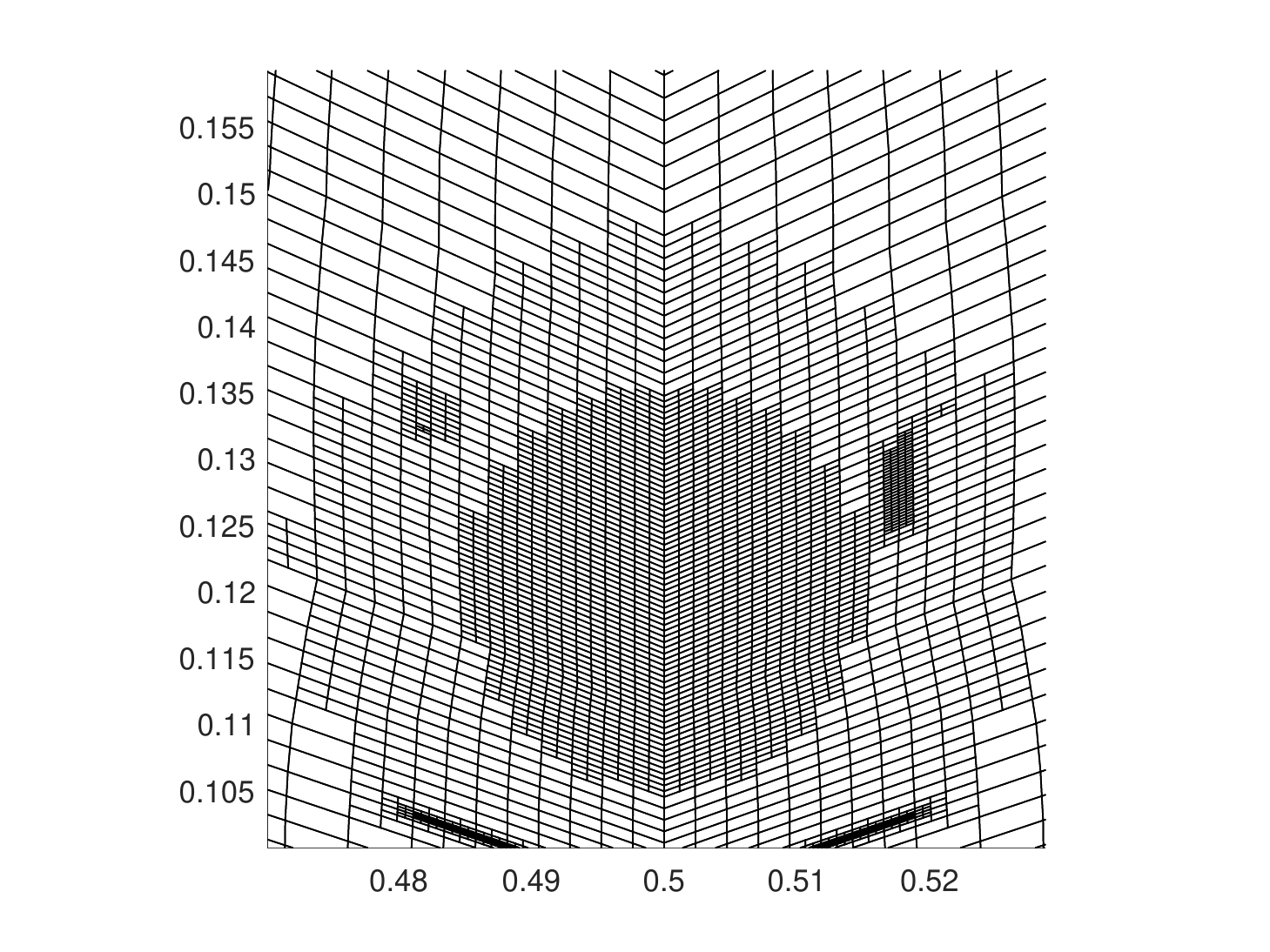}}\\
  \caption{The initial and the final adapted meshes. (a) The initial mesh;
  (b) After 20 refinement steps; (c) After 30 refinement steps; (d) The zoomed mesh in (c)}\label{refineVEMmesh}
\end{figure}

\begin{figure}[!htb]
  \centering
  \subfigure[Exact]{\includegraphics[scale=0.45]{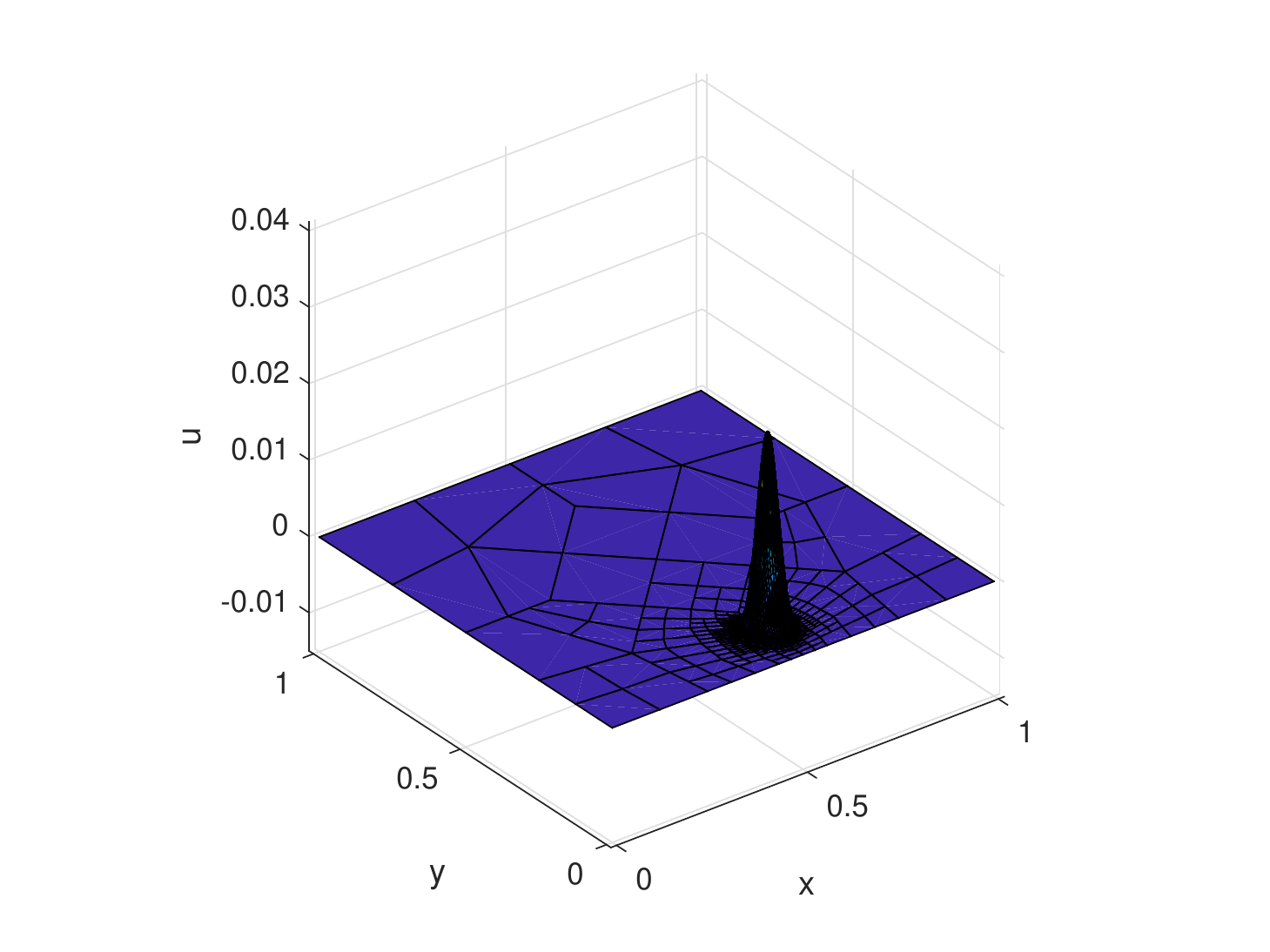}}
  \subfigure[Numerical]{\includegraphics[scale=0.45]{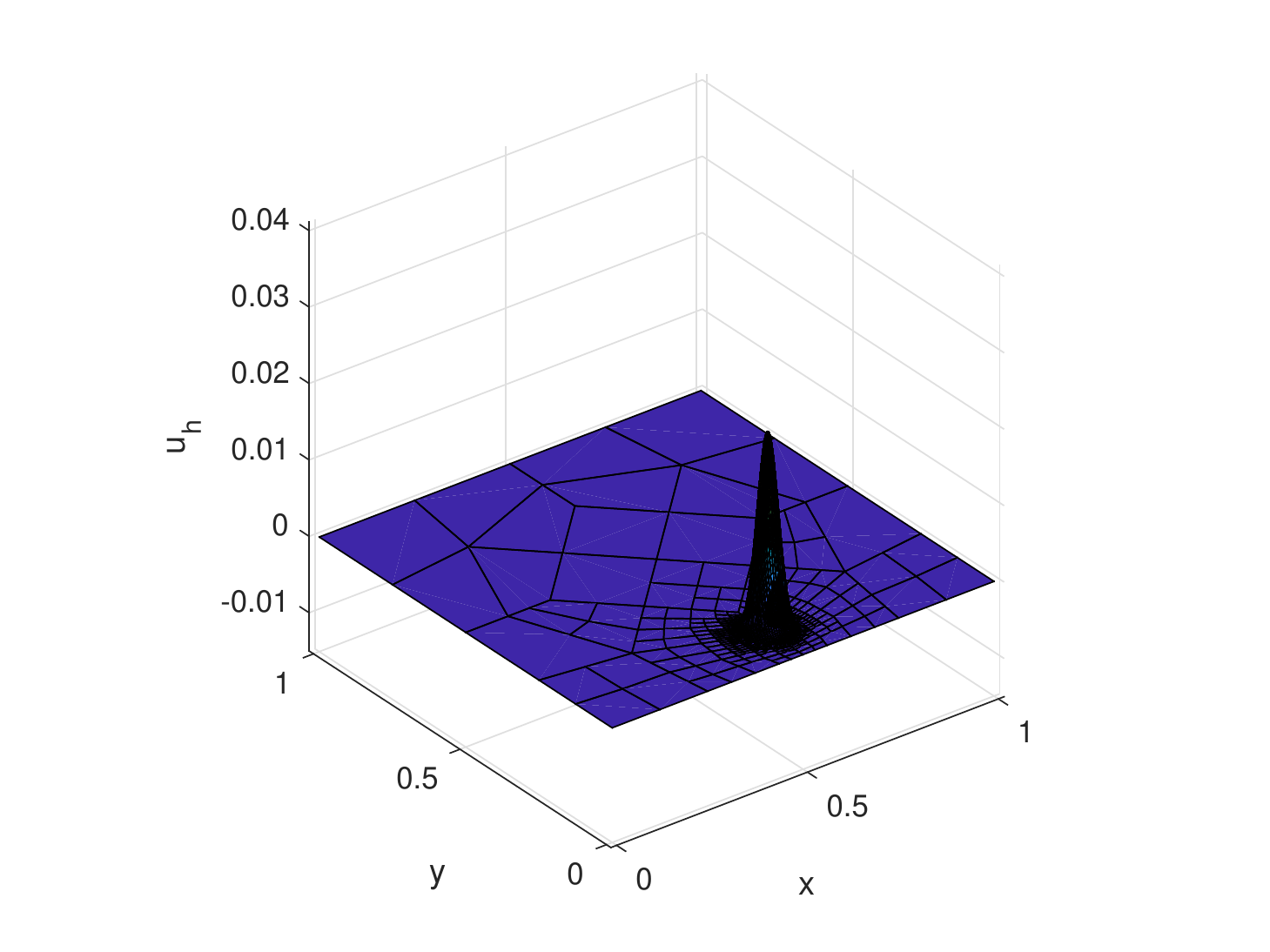}}\\
  \caption{The exact and numerical solutions}\label{refineVEMsol}
\end{figure}

The initial mesh and the final adapted meshes after 20 and 30 refinement steps are presented in Fig.~\ref{refineVEMmesh}~(a-c),
respectively. The detail of the last mesh is shown in Fig.~\ref{refineVEMmesh}~(d). Clearly, no small edges are observed. We also plot the adaptive approximation in Fig.~\ref{refineVEMsol}, which almost coincides with the exact solution. The full code is still available from mVEM package. The subroutine \mcode{PoissonVEM\_indicator.m} is used to compute the local indicators and the test script is \mcode{main\_Poisson\_avem.m}.

\end{document}